\baselineskip=15pt plus 2pt
\magnification =1200
\baselineskip=15pt plus 2pt
\magnification =1200
\def\sqr#1#2{{\vcenter{\vbox{\hrule height.#2pt\hbox{\vrule width.#2pt
height#1pt\kern#1pt \vrule width.#2pt}\hrule height.#2pt}}}}
\def\square{\mathchoice\sqr64\sqr64\sqr{2.1}3\sqr{1.5}3}
 at 10truept
\font\smalletters=cmr8 at 10truept
 at 10truept
 at 10truept at 12truept

\font\medtenrm=cmr10 scaled\magstep2
\centerline {\medtenrm Hankel operators that commute with
second-order differential operators}\par
\vskip.05in
\centerline {\medtenrm Gordon Blower}\par
\vskip.05in
\centerline  {\sl Department of Mathematics and Statistics, 
Lancaster University}\par
\centerline  {\sl Lancaster, LA1 4YF, England, UK. E-mail: 
g.blower@lancaster.ac.uk}\par
\vskip.1in
\centerline {7th October 2007}\par
\vskip.1in
\hrule
\vskip.1in

\noindent {\bf Abstract}\par
\indent Suppose that $\Gamma$ is a continuous and
self-adjoint Hankel operator on $L^2(0, \infty )$ with kernel $\phi
(x+y)$ and that
$Lf=-{{d}\over{dx}}(a(x){{df}\over{dx}})+b(x)f(x)$ with $a(0)=0$. 
If $a$ and $b$ are both quadratic, hyperbolic or trigonometric
functions, and $\phi$ satisfies a suitable form of Gauss's hypergeometric differential 
equation, or the confluent hypergeometric equation, then
$\Gamma L=L\Gamma$. The
paper catalogues the commuting pairs $\Gamma$ and $L$, including
important cases in random matrix theory. There are also results
proving rapid decay of the singular numbers of Hankel integral 
operators with kernels that are analytic and of exponential decay in 
the right half-plane.\par
\vskip.05in
\noindent MSC2000 Classification 47B35\par
\vskip.05in
\vskip.05in
\hrule
\vskip.1in
\noindent {\bf 1. Introduction}\par
\vskip.05in
\indent For $\phi\in L^2(0, \infty )$, we recall that the Hankel
 operator $\Gamma_\phi$ with kernel $\phi (x+y)$ is the integral operator 
$$\Gamma_\phi f(s)=\int_0^\infty \phi (s+t)f(t)\, dt\eqno(1.1)$$
\noindent from a subspace of $L^2(0, \infty )$ into $L^2(0,\infty
)$.\par
\indent Megretski\u\i, Peller and Treil [6] determined the possible spectral
multiplicity function that a continuous and self-adjoint Hankel 
operator can have; however, their spectral theory does not yield much
information about the eigenvectors. Tracy and Widom observed that some
self-adjoint and compact Hankel operators commute with self-adjoint
second-order differential operators $L$ that have purely discrete spectrum, and hence $\Gamma_\phi$ and
$L$ have a common orthonormal basis of eigenfunctions; see [9, 10].
One can then use the WKB approximation to describe the
asymptotic eigenvalue distribution of $L$,  and the asymptotic behaviour of
the eigenfunctions.\par
\vskip.05in
\noindent ------------------\par
\vskip.05in
\noindent {\smalletters This work was partially supported by EU Network Grant MRTN-CT-2004-511953
`Phenomena in High Dimensions'.}\par
\vfill
\eject

\indent In this paper we start with $L$ and seek choices of $\phi$, thus
reversing the order of steps in [9, 10]. Suppose that 
$$Lf=-{{d}\over{dx}}\Bigl(a(x){{df}\over{dx}}\Bigr)+b(x)f(x),\eqno(1.2)$$
\noindent where $a$ and $b$ are three times differentiable functions    
such that $a(0)=0$. Suppose that
 $a'''=Aa'$ and $b'''=Bb'$, for some constants $A$ and $B$. In section
2 we derive an explicit differential equation 
$$\phi''(u)+\alpha (u)\phi'(u)-\beta (u)\phi (u)=0\eqno(1.3)$$
\noindent which ensures that $\Gamma_\phi L=L\Gamma_\phi .$\par 
\indent The possibilities for $a$ and $b$
depend upon the sign of $A$, and in sections 3, 4 and 5 we consider the
quadratic, hyperbolic and trigonometric cases in detail. In all cases, the
differential equation reduces to a linear differential equation with
rational functions as coefficients which has less than or equal to
three singular points, and we determine the nature of the
singularities. Thus we prove the following theorem.
\vskip.05in
\noindent {\bf Theorem 1.1.} {\sl The differential equation (1.3) may be
transformed by change of variables to the hypergeometric equation
$$x(1-x){{d^2\phi}\over{dx^2}}+\{ \lambda -(\mu +\nu +1)x\}
{{d\phi}\over{dx}}-\mu\nu \phi =0\eqno(1.4)$$
or confluent
hypergeometric equation.}\par
\vskip.05in
\indent This result covers cases relating to standard models in
physics. In [11] Tracy and Widom considered the integral
operators that have kernels of the form
$$W(x,y)={{f(x)g(y)-f(y)g(x)}\over{x-y}}\eqno(1.5)$$
\noindent where $f$ and $g$ satisfy
$$m(x){{d}\over{dx}}\left[\matrix{f(x)\cr g(x)\cr}\right] =
\left[\matrix{ p(x)&q(x)\cr -r(x)&
-p(x)\cr}\right]\left[\matrix{f(x)\cr g(x)\cr}\right]\eqno(1.6)$$
\noindent for some real polynomial functions $p(x),q(x), r(x)$ 
and $m(x)$. Of particular interest are the cases where $m(x)=1, x$ or
$1-x^2$. The main application of Tracy--Widom operators $W$ is in random
matrix theory, where they describe the eigenvalue distributions of
$n\times n$ random matrices from the generalized unitary ensemble as
$n\rightarrow\infty$.\par
\indent Tracy and Widom introduced an indirect process for computing the
spectrum of $W$, and applied it effectively to the Airy and Bessel kernels
which are fundamentally important cases in random matrix theory; see
[9, 10].  Specifically, the edge distribution is given by
$\det (I-WP_{(x, \infty )})$, where $P_{(x, \infty )}$ is the
orthogonal projection $L^2(0, \infty )\rightarrow L^2(x, \infty )$ for
$x>0$. Their
first step was to introduce a self-adjoint Hankel operator $\Gamma$
such that $\Gamma^2=W,$ then to use the spectral theory of Hankel
operators to deduce information about the spectrum of $W$.\par
\indent The identity $W=\Gamma^2$ is equivalent to the factorization of
kernels 
$${{f(x)g(y)-f(y)g(x)}\over {x-y}}=
\int_0^\infty \phi (x+t)\phi(y+t)\, dt.\eqno(1.7)$$
\noindent Now $\Gamma_\phi$
is Hilbert--Schmidt if and only if $\int_0^\infty t\vert\phi
(t)\vert^2dt <\infty $. If $\Gamma_\phi$ is Hilbert--Schmidt, then
$\Gamma_\phi^2$ is of trace class, and $\det (I-\Gamma_\phi^2)$ is defined.\par 
\indent In [2, 3] we considered several differential equations of the
form (1.5) and
resolved whether factorization of $W$ takes place, giving explicit formulae
for $\phi$; several important examples involve the confluent
hypergeometric equation, which may be reduced by change of variable
to Whittaker's equation. For a more general statement about reducing 
systems, and how changes of variable can affect factorization, see
[3, section 3].\par
\indent By carrying out the reduction
in Theorem 1.1 explicitly, we recover the Tracy--Widom kernels as special cases of a 
more general theory, and also obtain new examples of commuting 
$\Gamma_\phi$ and $L$; 
there remain a few residual cases that we have not been able to solve
explicitly in terms of standard functions.\par
\indent In section 6 we present results concerning the
rate of decay of the singular numbers of $\Gamma$ under hypotheses
that are well suited to applications in random matrix theory.\par
\indent We follow standard notation as used in [4], and let $L^{(\alpha )}_n$ be the Laguerre polynomial of degree $n$
and order $\alpha$, so
$$L^{(\alpha
)}_n(s)={{s^{-\alpha}}\over{n!}}e^s{{d^n}\over{ds^n}}\bigl( e^{-s}
s^{n+\alpha }\bigr)\qquad (s>0).\eqno(1.8)$$
\vskip.05in
\noindent {\bf 2. The main result}\par
\vskip.05in
\noindent {\bf Theorem 2.1.} {\sl Let $\Gamma$ be the Hankel operator that has
kernel $\phi (x+y)$, let $L$ be the differential operator 
$$Lf(x)=-{{d}\over{dx}}\Bigl( a(x){{df}\over{dx}}\Bigr)
+b(x)f(x),\eqno(2.1)$$
\noindent where $a'''(x)=Aa'(x)$ and $a(0)=0$.\par
\indent (i) Then there exists
a real function $\alpha$ such that $\alpha
(x+y)=(a'(x)-a'(y))/(a(x)-a(y))$.\par
\indent (ii) Suppose further that $\beta
(x+y)=(b(x)-b(y))/(a(x)-a(y))$ for some real function $\beta$ and that 
$$\phi''(u)+\alpha (u)\phi' (u)-\beta (u)\phi (u)=0.\eqno(2.2)$$ 
\indent Then the operators $\Gamma$ and $L$ commute on $C_c^\infty
(0, \infty )$.}\par
\vskip.05in
\noindent {\bf Proof.} Let $T_t:L^2(0, \infty )\rightarrow L^2(0,
\infty )$ be the translation operator $T_tf(x)=f(x+t)$ and
$T_t^\dagger$ the adjoint for $t>0$. Then by the fundamental property
of Hankel operators, we have $T_t\Gamma =\Gamma T_t^\dagger $ and
hence ${{\partial }\over{\partial x}}\Gamma =-\Gamma{{\partial
}\over{\partial x}}$ on $C_c^\infty (0, \infty )$. To exploit this, we introduce the expression
$$\Phi (x,y)=(a(y)-a(x))\phi''(x+y)+(a'(y)-a'(x))
\phi'(x+y)-(b(y)-b(x))\phi(x+y);\eqno(2.3)$$
\noindent then by successive integrations by parts, we obtain
$$(L\Gamma -\Gamma L)f(x)=\Bigl[\phi
(x+y)a(y)f'(y)-\phi'(x+y)a(y)f(y)\Bigr]_0^\infty 
+\int_0^\infty \Phi (x,y)f(y)\, dy.\eqno(2.4)$$
\noindent The term in square brackets vanishes since $a(0)=0$; so the
operators $\Gamma $  and $L$ commute if and only if 
$\Phi=0$. The main idea is to reduce the
 condition $\Phi =0$ to the differential
equation 
$$\bigl(a(x)-a(y)\bigr)\bigl(\phi''(x+y)+\alpha (x+y)\phi'(x+y)-\beta (x+y)\phi
(x+y)\bigr)=0,\eqno(2.5)$$
\noindent involving
functions $\alpha$ and $\beta$ of only one variable $u=x+y$. The following lemma guarantees the
existence of such $\alpha$ and $\beta$ under suitable hypotheses,
and hence gives the proof of Theorem 2.1.\par 
\vskip.05in

\noindent {\bf Lemma 2.2.} {\sl Suppose that $a$ and $b$ are three times
continuously differentiable and non constant real functions.\par
\noindent (i) There exists a differentiable function $\alpha$ such that 
$$\alpha (x+y)={{a'(x)-a'(y)}\over
{a(x)-a(y)}}\eqno(2.6)$$
\noindent if and only if $a'''(x)=Aa'(x)$ for some constant $A$.\par
\noindent (ii) If $a'''(x)=Aa'(x)$, then $a(x)-a(y)=f(x+y)g(x-y)$ for
some differentiable functions $f$ and $g$.\par
\noindent (iii) Suppose that, for some differentiable functions $g$ and $h$, 
$${{b(x)-b(y)}\over{g(x-y)}}=h(x+y).\eqno(2.7)$$
\noindent  Then there exists a constant $B$ such that 
$b'''(x)=Bb'(x)$.}\par

\vskip.05in
\noindent {\bf Proof.} (i) $(\Rightarrow)$ We have 
$$\eqalignno{0=&\Bigl( {{\partial }\over{\partial x}}-
{{\partial }\over{\partial y}}\Bigr)\Bigl({{a'(x)-a'(y)}\over
{a(x)-a(y)}}\Bigr)&(2.8)\cr
&={{(a''(x)a(x)-a'(x)^2)-(a''(y)a(y)-a'(y)^2)+(a(x)a''(y)-a''(x)a(y))}
\over{(a(x)-a(y))^2}},&(2.9)\cr}$$
\noindent so we have
$${{\partial^2}\over{\partial x\partial y}}\Bigl(
a(x)a''(y)-a''(x)a(y)\Bigr)=0,\eqno(2.10)$$
\noindent which soon reduces to $a'''(x)=Aa'(x)$ for some constant
$A$.\par 
\indent (i) $(\Leftarrow )$ There are three main families of solutions of the
differential equation $a'''(x)=Aa'(x)$, depending upon the sign of
$A$:\par

\indent (Q) quadratic case when $A=0$ and $a(x)=q_2x^2+q_1x$, with
$q_2$ and $q_1$ not both zero, then 
$$\alpha (u)={{2q_2}\over{q_2u+q_1}};\eqno( 2.11)$$
\indent (H) hyperbolic case where $A>0$ and $a(x)=h_1\cosh tx+h_2\sinh
tx+h_3$, with $t^2=A$, $h_3=-h_1$ and $h_1$ and $h_2$ not both zero, 
then $a(0)=0$ and
$$\alpha (u)={{th_1\cosh tu/2+th_2\sinh tu/2}\over {h_1\sinh
tu/2+h_2\cosh tu/2}};\eqno(2.12)$$
\indent (C) circular case when $A<0$ and $a(x)=c_1\cos tx+c_2\sin tx+c_3$, 
where $t^2=-A$, $c_3=-c_1$, and $c_1$ and $c_2$ not both zero; then
$a(0)=0$ and 
$$\alpha (u)={{-tc_1\cos tu/2-tc_2\sin tu/2}\over{-c_1\sin tu/2+c_2\cos
tu/2}}.\eqno(2.13)$$
\indent (ii) One considers cases (Q),(H) and (C), and applies the
addition rule of trigonometry or hyperbolic trigonometry to
factorize $a(x)-a(y)$.\par 
\indent (iii) When such $g$ and $h$ exist, we have
$$\Bigl( {{\partial }\over{\partial x}}-
{{\partial }\over{\partial y}}\Bigr)\Bigl({{b(x)-b(y)}\over
{g(x-y)}}\Bigr)=0,\eqno(2.14)$$
\noindent which reduces to 
$${{b'(x)+b'(y)}\over
{b(x)-b(y)}}={{2g'(x-y)}\over{g(x-y)}},\eqno(2.15)$$
\noindent and arguing as in (i), we deduce that $b'''(x)=Bb'(x)$ for
some constant $B$.\par
\rightline{$\square$}\par
\vskip.1in
\noindent {\bf Remark.} Without loss of generality we can assume that $a$ and $b$ are
bounded above or below on $(0, \infty )$. If $a(x)\geq 0$ for all $x>0$
and $b$ is bounded below, then the quadratic form associated with $L$
is bounded below on $C_c^\infty (0, \infty )$, and hence $L$ admits of
a self-adjoint extension by Friedrichs's theorem.\par
\vskip.05in

\noindent {\bf 3. Quadratic cases}\par
\indent In this section, we determine all
the quadratic and linear choices of $a$ and $b$ and the corresponding $\phi$
such that $L$ and $\Gamma_\phi$ commute. We look for $\alpha$ and
$\beta$ as in Theorem 2.1.\par
\vskip.05in
\noindent {\bf Lemma 3.1.} {\sl Suppose that $a(x)=q_2x^2+q_1x+q_0$. Then
$\beta$ of Theorem 2.1 exists if and only if $b$ is quadratic or
linear; so that 
$b(x)=b_2x^2+b_1x+b_0$, and then}
$$\beta (x)={{b_2x+b_1}\over{q_2x+q_1}}.\eqno(3.1)$$
\vskip.05in
\noindent {\bf Proof.} We need to find
$\beta$ and $b$ such that 
$${{b(x)-b(y)}\over{x-y}}=(q_2(x+y)+q_1)\beta (x+y),\eqno(3.2)$$
\noindent so $b'''(x)=Bb'(x)$ by Lemma 2.2. Hence $b$ must belong to one of
the types (Q), (H) and (C), and evidently by (2.15) only a quadratic has the
right form.\par
\rightline{$\square$}\par
\vskip.05in
\noindent {\bf Proposition 3.2.} {\sl The differential equation (2.2) for
$\phi$ is}
$$-\phi''(u)-{{2q_2}\over{q_2u+q_1}}\phi'(u)+{{b_2u+b_1}\over{q_2u+q_1}}
\phi (u)=0,\eqno(3.3)$$
\noindent {\sl which may be transformed by change of variables to an elementary or confluent
form of the hypergeometric differential equation.}\par
\vskip.05in
\noindent {\bf Proof.} This formula follows from Theorem 2.1 and 
Lemma 3.1. The  
differential equation (3.3) has less than or equal to two singular
points, the only possibilities being $-q_1/q_2$ and infinity. Determining the
nature of the singularities in all possible cases below, we find that
 cases Q(i) and Q(ii) are trivial; cases Q(iii) and Q(iv)
are elementary; in cases
Q(iv)-(viii), infinity is an irregular singular point; in cases
Q(vi)-(viii) $-q_1/q_2$ is a regular singular point. We deduce that as
 a special
case of [12, p. 352] and 
[4, p. 1084], the solution is either an elementary function or may be expressed in terms of confluent
hypergeometric functions. \par
\rightline{$\square$}\par
\vskip.05in
\noindent {\bf Examples 3.3.} {\sl (Quadratic cases).} We now carry out a systematic reduction of the
(3.3), bearing in mind that we wish to have a continuous Hankel
operator $\Gamma_\phi$.\par 
\indent Q(i) If $q_1=q_2=0$, then $a(x)=0$ and $b(x)$ is constant by
(2.3).\par
\indent After translating the variable $u$, we assume henceforth 
without losing generality that exactly one of $q_1$ and $q_2$ is non-zero.\par 
\indent Q(ii) If $q_2=b_2=b_1=0$, then $\phi$ is linear, so does not
give a continuous Hankel operator.\par
\indent Q(iii) If $b_2=b_1=0$, $q_2>0$ and $q_1\geq 0$, then 
$$\phi (u)=C_1+{{C_2}\over{q_2u+q_1}}.\eqno(3.4)$$ 
\indent In particular, when $q_2=1$, $q_1=0$, $C_1=0$ and $C_2=1$,  
we thus obtain Carleman's operator $\Gamma$ with kernel $1/(x+y)$. Carleman's operator is continuous on $L^2(0, \infty )$ and has
spectrum $[0, \pi ]$ with multiplicity two by a result due to Power [6].
 Both $\Gamma$ and the
differential operator $L=-{{d}\over{dx}}\bigl(x^2{{d}\over{dx}}\bigr)$ are multipliers for
the Mellin transform, so this case of Theorem 2.1 was to be
expected.\par
\indent Q(iv) If $q_2=b_2=0$ and $b_1, q_1\neq 0$, then $\phi$ is
hyperbolic or trigonometric according to the sign of $b_1/q_1$. In
particular, $b_1=q_1$ gives $\phi (u)=e^{-u}$, and $\Gamma_\phi$ is
a rank-one continuous Hankel operator. We defer the case of
$b_1/q_1<0$ until section 5.\par
\indent Q(v) If $q_2=0$ and $b_2\neq 0$, then we obtain the Airy
equation. In particular, if $a(x)=x$ and $b(x)=x(x+s)$, then $\phi$ 
satisfies $\phi''(u)=(u+s)\phi (u),$ so
 the Airy function $\phi (u)={\hbox{Ai}}(u+s)$ gives a solution. Here
 the associated Hankel operator $\Gamma_\phi$ is Hilbert--Schmidt, and 
 $L$ has discrete spectrum as described by Titchmarsh's theory of
oscillations [8]; Tracy and Widom discuss the asymptotics of the
common sequence of eigenfunctions in various ranges [9]. The Hankel
operator $\Gamma_\phi$ has square 
$${{{\hbox{Ai}}(x){\hbox{Ai}}'(y)-{\hbox{Ai}}'(x){\hbox{Ai}}(y)}\over{x
-y}},\eqno(3.5)$$
\noindent which is the famous Airy kernel as used in [9] to describe the soft
edge of eigenvalues distributions from the Gaussian unitary
ensemble.\par
\indent When $q_1=1$ and $b_2<0$,
we obtain $\phi''(u) +u\phi (u)=0$ by translation and rescaling, and 
the general solution of this is
$$\phi (u)=C_1u^{1/2}J_{1/3}\Bigl({{2}\over{3}}u^{3/2}\Bigr) +
C_2u^{1/2}J_{-1/3}\Bigl({{2}\over{3}}u^{3/2}\Bigr),\eqno(3.6)$$
\noindent for typical $C_1$ and $C_2$ with asymptotics [5, p.133]
$$\phi (u)\asymp Cu^{-1/4}\cos\bigl( {{2}\over{3}}u^{3/2}+c\bigr)\qquad
(u\rightarrow\infty ).$$
 \indent Q(vi) If $b_2=0$ and $b_1,q_2\neq 0$, then we obtain after
translation and rescaling the equation
$$-\phi''(u)-{{2}\over{u}}\phi' (u)\pm {{1}\over{4u}}\phi
(u)=0\eqno(3.7)$$
\noindent with solutions including [5, p.250]
$$\phi_+(u)=u^{-1/2}K_1(\sqrt {u}),\qquad \phi_-(u)=u^{-1/2}J_1(\sqrt
{u}).\eqno(3.8)$$
\noindent We have by [5, p. 435]
$$\phi_+(u)\asymp\cases{ {{1}\over{u}}& as $u\rightarrow 0+$\cr
                         \Bigl({{\pi}\over{2}}\Bigr)^{1/2}
u^{-3/4}e^{-\sqrt{u}}& as $u\rightarrow\infty$,\cr}\eqno(3.9)$$
\noindent so $\Gamma_{\phi_+}$ defines a continuous linear operator on
$L^2(0, \infty )$ which is not Hilbert--Schmidt; but the compression of
$\Gamma_{\phi_+}$ to $L^2(x, \infty )$ is Hilbert--Schmidt for $x>0$. 
Further by [5],
$$\phi_-(u)\asymp\cases{ {{1}\over{2}}& as $u\rightarrow 0+$\cr
                         \Bigl({{\pi}\over{2}}\Bigr)^{1/2}
u^{-3/4}\cos (\sqrt{u}-3\pi /4)& as
$u\rightarrow\infty$,\cr}\eqno(3.10)$$
\noindent so $\Gamma_{\phi_-}$ is not Hilbert--Schmidt. As in [3,
2.3], the
associated integrable operator is
$${{x\phi_- (x)(y\phi_- (y))'-(x\phi_- (x))'y\phi_-
(y)}\over{x-y}}={{1}\over{4}}\int_0^\infty \phi_- (x+t)\phi_- (y+t)\, dt.$$

\indent Q(vii) If $b_2,q_2\neq 0$, then we obtain after
translation the equation 
$$-\phi''(u)-{{2}\over{u}}\phi' (u)+\Bigl(
{{b_2}\over{q_2}}+{{b_1}\over{q_2u}}\Bigr)\phi (u)=0.\eqno(3.11)$$
\noindent In particular, when $b_1=0$, we can rescale to the equation
$$-\phi''(u)-{{2}\over{u}}\phi' (u)\pm \phi (u)=0$$
\noindent with solutions
$$\phi_+(u)={{e^{-u}}\over{u}},\qquad \phi_-(u)={{\kappa_1\cos u+\kappa_2\sin u}\over {u}}.\eqno(3.12)$$
\noindent By Theorem 2.1 of [3], with $f(x)=e^{-x}$ and 
$$g(x)=e^{-x}\int_0^\infty {{e^{-2t}dt}\over{x+t}},$$
\noindent we have
$${{f(x)g(y)-f(y)g(x)}\over{x-y}}=\int_0^\infty
\phi_+(x+t)\phi_+(y+t)\, dt.\eqno(3.13)$$
\noindent This $\Gamma_{\phi_+}$ defines a continuous linear operator on
$L^2(0, \infty )$, such that the compression to $L^2(x, \infty )$ is
Hilbert--Schmidt for $x>0$. \par
\indent Further, $\Gamma_{\phi_-}$ defines a continuous
linear operator on $L^2(0, \infty )$ since the kernel is a sum of Schur
multiples of the Carleman operator, namely
$$\phi_-(x+y)=\kappa_1\Bigl({{\cos x\cos
y}\over {x+y}}- {{\sin x\sin y}\over {x+y}}\Bigr)+
\kappa_2\Bigl({{\sin x\cos
y}\over {x+y}}+ {{\cos x\sin y}\over {x+y}}\Bigr).\eqno(3.14)$$ 
\indent With $f(x)=e^{ix}$ and 
$$g(x)=e^{ix}\int_0^\infty {{e^{2it}dt}\over{x+t}},\eqno(3.15)$$
\noindent we have by Theorem 2.1 of [3] the identity
$${{f(x)g(y)-f(y)g(x)}\over{x-y}}=\int_0^\infty
{{e^{i(x+y+2t)}dt}\over{(x+t)(y+t)}}.\eqno(3.16)$$

\indent Q(viii) Now we consider (3.11) when $q_2=b_2=1$, and
$b_1=-2(n+1)$. Then $g(u)=u\phi (2^{-1}u)$ satisfies Laguerre's
equation
$$g''(u)+\Bigl({{-1}\over{4}} +{{n+1}\over{u}}\Bigr)
g(u)=0,\eqno(3.17)$$
\noindent so that, when $n+1$ is a positive integer, $g(u)=ue^{-u/2}L^{(1)}_n(u)$ and $\phi
(u)=e^{-u}L_{n}^{(1)}(2u)$ gives a solution. See 2.2 of [3] for the
corresponding Tracy--Widom operator.\par
\indent This gives a complete catalogue of the possible quadratic
cases.\par
\vskip.05in
\noindent {\bf 4. Hyperbolic cases}\par
\indent In this section we consider the case in which $a$ and $b$ are
hyperbolic functions, giving some $L$, and obtain
$\Gamma$ in terms of standard special functions such that $L$ and
$\Gamma$ commute.\par
\indent Without loss of generality we can choose $t=2$, since other 
cases occur by rescaling. The change of variables 
$x=e^{-2u}$ gives a unitary transformation $L^2((0, \infty );
du)\rightarrow L^2((0,1); dx/x)$, and we modify the
definition of the Hankel operator accordingly.\par
\vskip.05in

\noindent {\bf Definition} {\sl (Hankel operator).} For
$\rho\in L^2((0,1);dx/x)$, the Hankel operator 
with kernel $\rho (xy)$ is
$$\Gamma_\rho h(x)=\int_0^1\rho (xy)h(y){{dy}\over{y}},\eqno(4.1)$$
\noindent where $h$ is in some subspace of $L^2((0,1); dy/y)$.\par
\vskip.05in
\vskip.05in
\noindent {\bf Lemma 4.1.} {\sl Suppose that $a$ is
hyperbolic, so $a(x)=h_1\cosh tx
+h_2\sinh tx +h_3$ where $h_1$ and $h_2$ are not both zero. Then
$\beta$ of Theorem 2.1 exists if and only if $b(x)=h_4\cosh tx+h_5\sinh tx+h_6$ for
some constants and then}
$$\beta (v)= {{h_4\sinh tv/2+h_5\cosh tv/2}\over{h_1\sinh tv/2+h_2\cosh
tv/2}}.\eqno(4.2)$$
\vskip.05in
\noindent {\bf Proof.} This follows from Lemma 2.2 since we need to find
$\beta$ and $b$ such that 
$${{b(x)-b(y)}\over{\sinh t(x-y)/2}}=2\bigl(h_1\sinh t(x+y)/2+h_2\cosh
t(x+y)/2\bigr)\beta (x+y).\eqno(4.3)$$
\noindent Hence $b'''(x)=Bb'(x),$ and of the types (Q), (C) and (H),
only a hyperbolic $b$ has the right form with $v=x+y$.\par
\rightline{$\square$}\par
\vskip.05in

\vskip.05in
\noindent {\bf Proposition 4.2.} {\sl (i) With $t=2$ and $x=e^{-2u}$, 
the differential equation (2.2) for
$\phi$ becomes in the new variable}
$${{d^2\phi}\over{dx^2}}+{{2(h_2-h_1)}\over{h_1+h_2+(h_2-h_1)x}}
{{d\phi}\over{dx}}-{{h_4+h_5+(h_5
-h_4)x}\over {4x^2(h_1+h_2+(h_2-h_1)x)} }\phi =0,\eqno(4.4)$$
\noindent {\sl and the commuting differential operator $L$ becomes}
$$-{{d}\over{dx}}\Bigl(2x\bigl\{
(h_1-h_2)x^2+2h_3x+(h_1+h_2)\bigr\}{{df}\over{d
x}}\Bigr)+\Bigl({{x}\over{2}}(h_4-h_5)+h_6+{{1}\over{2x}}(h_4+h_5)\Bigr)
f.$$
\noindent {\sl (ii) The equation (4.4) may be reduced by change of
variables to the hypergeometric equation or
the confluent hypergeometric equation.}\par 
\vskip.05in
\noindent {\bf Proof.} (i) We have  
$$\alpha
={{2h_1(1+x)+2h_2(1-x)}\over{h_1(1-x)+h_2(1+x)}}\qquad
(0<x<1),\eqno(4.5)$$
\noindent and by Lemma 4.1 we have
$$\beta
={{h_4(1-x)+h_5(1+x)}\over{h_1(1-x)+h_2(1+x)}}\qquad
(0<x<1),\eqno(4.6)$$
\noindent so (2.2) reduces as stated. One obtains the formula for 
$L$ by changing variables.\par
\indent (ii) In 4.3 below, we consider the nature of the singular
points of the differential equation. Effectively there are four
constants in (4.4), namely $h_1\pm
h_2$ and $h_4\pm h_5$; by taking
$h_3=-h_1$, we ensure that $a(0)=0$. One can easily verify that the
effect of the change of variable $x=1/y$ is to preserve the shape of
the formula (4.4) and to interchange the constants 
$h_1+h_2\leftrightarrow h_2-h_1$
and $h_4+h_5\leftrightarrow h_5-h_4$.\par
\indent In cases H(i), H(iii) and H(vii) below, there are two regular 
singular points, so the solution is 
elementary; in cases H(ii) and H(iv),
the singular points are zero and infinity, and one of them is irregular, 
so the equation is confluent
hypergeometric type; whereas in the remaining cases H(v) and H(vi), 
the three singular points are all regular, so the equation is
of hypergeometric type.\par
\rightline{$\square$}\par
\vskip.05in

\noindent {\bf Examples 4.3} {\sl (Hyperbolic cases).} \par
\indent H(i) If $h_1=h_2$ and $h_4=h_5$, then the differential equation
(4.4) reduces to 
$$\phi''(u)-{{h_4}\over{4h_1u^2}}\phi (u)=0,\eqno(4.7)$$
\noindent which has solutions $\phi (u)=u^p$ where
$4p(p-1)=h_4/h_1.$ When $p>0$, the corresponding Hankel operator is
Hilbert--Schmidt.\par
\indent H(ii) If $h_1=h_2$ and $h_4\neq h_5$, we obtain after
rescaling
$$\phi''(u)+{{1-\nu^2}\over{4u^2}}\phi
(u)=\pm {{1}\over{u}}\phi (u),\eqno(4.8)$$
\noindent with $\nu >0$ where the solutions are 
$$\phi_+(u)=\sqrt {u}K_\nu (2\sqrt{u}),\qquad
\phi_-(u)=\sqrt{u}J_\nu(2\sqrt{u}).\eqno(4.9)$$  
\noindent Indeed, $\phi_+$ emerges for the choice $h_1=h_2=-1$, $h_3=1,$
$h_4=5-\nu^2$ and $h_5=-3-\nu^2$; whereas $\phi_-$ emerges
for $h_1=h_2=-1$, $h_3=1$, $h_4=-3-\nu^2$ and $h_5=5-\nu^2$. We have
$$\phi_+(u)\asymp\cases{
2^{-1}\Gamma (\nu )u^{(1-\nu )/2}& as $u\rightarrow 0+$\cr
                        2^{-1}\pi^{1/2}u^{1/4}e^{-2\sqrt{u}}& 
as $u\rightarrow\infty$,\cr}\eqno(4.10)$$
\noindent so $\Gamma_{\phi_+}$ defines a Hilbert--Schmidt operator on
$L^2((0, 1); dx/x)$ when $\nu <1$. \par
\indent Further, by [5, p.436] we have
$$\phi_-(u)\asymp\cases{\Gamma (\nu +1)^{-1}u^{(1+\nu )/2}& 
as $u\rightarrow 0+$\cr
                        \pi^{-1/2}u^{1/4}\cos\bigl( 2\sqrt{u}
-{{\pi\nu}\over{2}}-{{\pi}\over{4}}\bigr)& as
$u\rightarrow\infty$;\cr}\eqno(4.11)$$
\noindent so $\Gamma_{\phi_-}$ defines a Hilbert--Schmidt operator on
$L^2((0, 1); dx/x)$. The associated Tracy--Widom operator on $L^2((0,1); dy/y)$ has
kernel 
$${{\sqrt{x}J'_\nu (2\sqrt{x})J_\nu (2\sqrt{y})-J_\nu (2\sqrt{x}) 
\sqrt{y}J'_\nu (2\sqrt{y})}\over {x-y}}=\int_0^1 J_\nu(2\sqrt{tx}) 
J_\nu (2\sqrt{ty})dt,\eqno(4.12)$$
\noindent and the commuting differential operator is
$$L_\nu f(x)=4{{d}\over{dx}}\Bigl( x(1-x){{df}\over{dx}}\Bigr)
+\Bigl(-4x+{{1-\nu^2}\over {x}}-2\mu \Bigr) f(x)\eqno(4.13)$$
\noindent with boundary conditions $f(0)=f(1)=0.$ In random matrix
theory, this is associated with hard edges, such as occur with the
Jacobi ensemble.\par
\indent H(iii) If $h_1=-h_2$ and $h_5=-h_4$, then we obtain
$$\phi''(u)+{{2}\over{u}}\phi'(u)-{{h_5}\over{4h_2u^2}}\phi (u)=0,
\eqno(4.14)$$
\noindent so we have solutions $\phi (u)=u^p$ where
$4p^2+4p-h_5/h_2=0$.\par
\indent H(iv) If $h_1=-h_2$ and $h_4\neq -h_5$, then we change variables
to $\phi (x)=\psi (y)$ where $y=1/x$ and obtain
$${{d^2\psi}\over{dy^2}}+\Bigl({{h_4+h_5}\over{8h_1y}}+{{h_5-h_4}\over{
h_1y^2}}\Bigr)\psi (y)=0.\eqno(4.15)$$
\noindent As in case H(ii), we thus we obtain solutions
$$\phi_+(x)={{1}\over{\sqrt{x}}}K_\nu
\Bigl({{2}\over{\sqrt{x}}}\Bigr),\quad
\phi_-(x)={{1}\over{\sqrt{x}}}J_\nu \Bigl({{2}\over{\sqrt{x}}}\Bigr).$$
\noindent By (4.10), $\phi_+$ gives a Hilbert--Schmidt Hankel operator
on $L^2((0,1); dx/x);$ whereas, by (4.11), $\phi_-$ gives a 
Hankel operator which is
not Hilbert--Schmidt. \par
\indent H(v) When $h_1\neq 0$, $h_2=0$ and $h_4=-h_5\neq 0$, the
equation (4.4) reduces to hypergeometric equation.
In particular, if $h_1=1/4$, $h_2=0$, $h_4=-h_5$, and
$\mu$ and $\nu$ satisfy $\mu +\nu =1$ and $\mu\nu
=-2h_4$, then we have (1.4) with 
$\lambda =0$ and nonzero parameters $\mu $ and $\nu$; however, the usual series for
the hypergeometric function $F(\mu, \nu ,\lambda , x)$ is then undefined. By an identity from [4, p.
1073], the function
$$\phi (x)=xF(\mu+1, \nu +1, 2,x),\eqno(4.16)$$
\noindent gives a power series solution which is analytic for 
$\vert x\vert <1$, with $\phi (0)=0$, and diverges everywhere on the circle  $\{z: \vert z\vert =1\}$ by
[4, p. 1066].\par
\indent H(vi) In the generic case $h_1\neq \pm h_2$ and $h_4\neq \pm
h_5$, we introduce the regular singular point 
$\zeta =-(h_1+h_2)/(h_2-h_1)$ and parameters $\alpha_1,\alpha_2,\beta_1$
and $\beta_2$ by the
simultaneous quadratic equations
$$\alpha_1\alpha_2=-{{h_4+h_5}\over{4(h_1+h_2)}},
\qquad \beta_1\beta_2=-{{h_5-h_4}\over{4(h_2-h_1)}},$$
$$\alpha_1+\alpha_2=1, \qquad \beta_1+\beta_2=1.\eqno(4.17)$$
\noindent Then the differential equation (4.4) is given in Riemann's
notation [12, 5] by 
$$\phi (x)=P\left\{\matrix{ 0&\infty & \zeta &{}\cr
 \alpha_1&\beta_1& -1; &x\cr
\alpha_2&\beta_2&0&{}\cr}\right\},\eqno(4.18)$$
\noindent and hence by [5, p 156], may be reduced to 
the hypergeometric equation (1.4).\par  
\indent H(vii) Finally, when $h_1\neq h_2$ and $h_4=h_5=0$, as in
Q(iii) we have elementary
solutions 
$$\phi (x)=C_1+{{C_2}\over{(h_2-h_1)x+h_2+h_1}}.\eqno(4.19)$$ 
\vskip.1in

\noindent {\bf 5. Circular cases}\par
\vskip.05in
\indent In this section we consider the remaining case, namely when
$\phi$, $a$ and $b$ are circular functions. The results of this section
are quite analogous to those of H(v) and H(vi); although they are somewhat
contrived, since the notion of a Hankel integral
operator over the circle is not in common use.\par
\indent Suppose that $a$, $b$ and $\phi$ have period $2\pi /t$, and that
$a(0)=a(2\pi /t)=0$. We extend $f\in C_c^\infty (0, \infty )$ by
$f(x)=0$ for $x<0$, and introduce $F(y)=\sum_{k=-\infty}^\infty
f(y+2\pi k/t)$, which is $2\pi /t$ periodic. Evidently
$$\int_0^\infty \phi (x+y)f(y)dy=\int_0^{2\pi /t}\phi
(x+y)F(y)dy\eqno(5.1)$$
\noindent is also $2\pi /t$ periodic.\par
\vskip.05in
\noindent {\bf Definition} {\sl (Hankel operator).} The Hankel operator on $L^2((0, 2\pi /t);
dy)$ with kernel $\phi (x+y)$ 
 is
$$\Gamma F(x)=\int_0^{2\pi /t} \phi (x+y) F(y)\, dy.$$
\indent We consider $\Gamma$ as an
operator on $C^\infty (0, 2\pi /t)$, and look for a second-order
differential operator $L$ as in (1.2) such that $L\Gamma =\Gamma L$.  \par
\vskip.05in

\noindent {\bf Lemma 5.1.} {\sl Suppose that $a(x)=c_1\cos tx
+c_2\sin tx +c_3$ is circular, where $c_1$ and $c_2$ are not both zero. 
Then $\beta$ of Theorem 2.1 exists if and only if $b(x)=c_4\cos tx+c_5\sin tx+c_6$ for
some constants and then}
$$\beta (u)= {{-c_4\sin tu/2+c_5\cos tu/2}\over{-c_1\sin tu/2+c_2\cos
tu/2}}.\eqno(5.2)$$
\vskip.05in
\noindent {\bf Proof.} This follows from Lemma 2.2 since we need to find
$\beta$ and $b$ such that 
$${{b(x)-b(y)}\over{\sin t(x-y)/2}}=2\Bigl(-c_1\sin t(x+y)/2+c_2\cos
t(x+y)/2\Bigr)\beta (x+y).\eqno(5.3)$$
\noindent Hence $b'''(x)=Bb'(x),$ and of the types (Q), (H) and (C),
only a circular $b$ has the right form.\par
\rightline{$\square$}\par
\vskip.05in
\noindent {\bf Proposition 5.2.} {\sl (i) Let $\tau =\tan u$ and
$t=2$. Then differential equation (2.2) for
$\phi$ becomes in the new variable} 
$${{d^2\phi }\over{d\tau^2}}+ {{2c_1}\over{c_1\tau
-c_2}}{{d\phi}\over{d\tau}}-{{-c_4\tau +c_5}\over{-c_1\tau +c_2}}{{\phi
(\tau )}\over{(1+\tau^2)^2}}=0,\eqno(5.4)$$ 
\noindent {\sl and the commuting differential operator transforms to}
$$-(1+\tau^2){{d}\over{d\tau}}
\Bigl( \bigl\{c_1(1-\tau^2)+2c_2\tau
+c_3(1+\tau^2)\bigr\} {{df}\over{d\tau}}\Bigr)
+\Bigl(c_4{{1-\tau^2}\over{1+\tau^2}}+c_5{{2\tau}\over{1+\tau^2}}+c_6
\Bigr)f.$$
\noindent {\sl (ii) The equation (5.4) may be reduced by change of
variables to the hypergeometric
equation.}\par
\vskip.05in
\noindent {\bf Proof.} (i) This follows from Theorem 2.1 and Lemma 5.1 
by calculation.\par
\indent (ii) The differential equation has regular
singular points: at $c_2/c_1$, when $c_1\neq 0$; at $\infty$, when
$c_1=0$; and at $\pm i$, when $c_4\neq 0$ or $c_5\neq 0$. When
$c_4=c_5=0$, the equation has elementary solutions, as in H(vii) and Q(iii). The effect of the change of
variable $\tau =1/s$ is to preserve the shape of (5.4) and to
interchange the constants $c_1\leftrightarrow c_2$ and
$c_4\leftrightarrow c_5$. By [5, p. 156], the differential equation
may be reduced by change of variables to Gauss's hypergeometric
equation. \par
\rightline{$\square$}\par
\vskip.05in
\noindent {\bf Examples 5.3} {\sl (Circular case).} We present solutions of the differential
equations in the special case where the singular points are $0, \pm
i$. Let $c_1=-1$, $c_2=0$, $c_4=3$ and $c_5=0$; so that, the 
differential equation in the original variables is 
$$\phi'' (x)+2\cot x\,\phi'(x)+3\phi (x)=0,\eqno(5.5)$$
\noindent which has general solution
$$\phi (x)=c_7\cos x-c_8({\hbox{cosec}} \,x-2\sin
x).\eqno(5.6)$$
\indent The Hankel operator on $L^2((0,\pi ); dx)$ with kernel $\cos
(x+y)$ has rank two and eigenfunctions $\cos x$ and $\sin x$. In the
new variable, the solution is the algebraic function
$$\phi (\tau )={{c_7\tau +c_8(\tau^2-1)}\over{\tau\sqrt{\tau^2+1}}}.$$  
\noindent {\bf Proof of Theorem 1.1.} The Propositions 3.2, 4.2 and 5.2
cover the three cases that together give the proof of Theorem
1.1.\par
\rightline{$\square$}\par
\vskip.05in
\vskip.05in
\noindent {\bf 6. Singular numbers of Hankel integral operators}\par
\indent Let $\Gamma_\phi$ be the Hankel operator on $L^2(0, \infty )$ with
kernel $\phi (x+y)$. In applications to random matrix theory, the ultimate aim of the analysis is to prove properties of
$\det (I-\Gamma_\phi^2)$. If $W=\Gamma_\phi^2$, and $W$ has eigenvalues
$\lambda_1\geq\lambda_2\geq \dots$ listed according to multiplicity, then
the singular numbers of $\Gamma_\phi$ are $s_j=\lambda_j^{1/2}$ for
$(j=1,2, \dots )$.  In this section we prove results which show that
if $\phi$ is analytic on a suitable domain and satisfies
various growth bounds, then $\Gamma_\phi$ is a trace-class operator and 
that its singular numbers are of
rapid decay. An important feature of the results is that if $\phi
(z)$ satisfies the hypotheses, then the translated
functions $\phi_s(z)=\phi (s+z)$ also
satisfy them for all $s>0$, up to modified constants. This reflects the fact that the
compression of $\Gamma_\phi^2$ to $L^2(s, \infty )$ is unitarily
equivalent to $\Gamma_{\phi_s}^2$ on $L^2(0, \infty )$.\par 
\vskip.05in
\noindent {\bf Definition} {\sl (Hankel matrix).} Suppose that
 $(\gamma_j)_{j=1}^\infty$ is a
sequence in $\ell^2({\bf Z}_+)$. Then the densely defined linear 
operator $\Gamma$ with matrix $[\gamma_{j+k-1}]$ with
respect to the standard orthonormal basis of $\ell^2({\bf Z}_+)$ is a Hankel
operator. \par
\vskip.05in

\noindent {\bf Theorem 6.1.} {\sl Let $\phi{}$ be an analytic 
function on the
half plane $\{ z=x+iy: x>-\delta \}$ for some $\delta >0$, and 
suppose that there
exists $\varepsilon >0$ and $M_\varepsilon (\delta )>0$ such that
$$\vert \phi (z)e^{\varepsilon z}\vert \leq M_\varepsilon (\delta
)\qquad (\Re z>-\delta ).\eqno(6.1)$$
\noindent (i) Then $\Gamma_\phi$ has singular numbers $s_j$ such that 
$$j^ps_j\rightarrow 0\qquad (j\rightarrow\infty )\eqno(6.2)$$
\noindent for all integers $p>0.$}\par
\noindent {\sl (ii) For all integers $p>1$, there exists $C_p>0$ such
that}
$$\log\det (I+x\Gamma_\phi ^2)\leq C_px^{1/p}\qquad (x>0),\eqno(6.3)$$
\noindent {\sl and the entire function $\det (I+z\Gamma_\phi^2)$ has order
zero.}\par
\vskip.05in
\noindent {\bf Proof.} (i) The operator with kernel 
$\varepsilon^{-1}\phi (\varepsilon^{-1}(x+y))$ is unitarily
equivalent to the operator with kernel $\phi (x+y)$, and
$\varepsilon^{-1}\phi (z/\varepsilon )$ satisfies (6.1) with
$\varepsilon =1$, $M_1(\varepsilon\delta )=\varepsilon^{-1}M_\varepsilon
(\delta )$ on $\{ z: \Re z>-\varepsilon \delta \}$. Hence we assume 
without loss of generality that $\phi$ has been rescaled so
that it satisfies (6.1) with
$\varepsilon =1$. Let $h_n(s)=e^{-s}L^{(0)}_n(2s)$, so that 
$(h_j)_{j=0}^\infty$
gives an orthonormal basis of $L^2(0, \infty )$, and is associated with
a unitary equivalence between $L^2(0, \infty )$ and $\ell^2({\bf
Z}_+)$. Then the Hankel
integral operator $\Gamma_\phi{}$ on $L^2(0, \infty )$ 
is unitarily equivalent to the Hankel
matrix $\Gamma = [\gamma_{j+k-1}]_{j,k=1}^\infty $ on $\ell^2({\bf
Z}_+)$, where $\gamma_j=\int_0^\infty
\phi{}(x)h_j(x)\, dx$ as in [6]. We shall show that the entries
of $\Gamma$ are of rapid decay as $j+k\rightarrow\infty$.\par
\indent By integrating repeatedly by parts, we find that
$$\eqalignno{\gamma_j&=\int_0^\infty \phi (x)h_j(x)\, dx\cr
  &={{(-1)^j}\over{j!}}\int_0^\infty \Bigl({{d^j}\over{dx^j}}
\bigl( e^x\phi (x)\bigr)\Bigr) x^je^{-2x}\, dx&(6.4)\cr}$$
\noindent where by (6.1) $e^z\phi (z)$ is analytic and bounded inside the circle
with centre $x$ and radius $x+\delta $. Applying Cauchy's estimates, we
obtain 
$$x^j\Bigl\vert {{d^j}\over{dx^j}}
\bigl( e^x\phi (x)\bigr)\Bigr\vert\leq j!M_1(\delta )\Bigl({{x}\over{x+\delta
}}\Bigr)^j\eqno(6.5)$$ 
\noindent and hence 

$$\vert \gamma_j\vert \leq M_1(\delta )\int_0^\infty {{x^j}\over{(x+\delta )^j}}
e^{-2x}\, dx,\eqno(6.6)$$
\noindent and summing this estimate over $j$ we obtain by the monotone
convergence theorem
$$\eqalignno{\sum_{j=1}^\infty j^{p+1}\vert \gamma_j\vert&\leq
M_1(\delta )\int_0^\infty\sum_{j=0}^\infty {{(j+p+1)!}\over{j!}}\Bigl(
{{x}\over{x+\delta }}\Bigr)^j e^{-2x}dx\cr
&\leq {{M_1(\delta )}\over{\delta^{p+2}}}\int_0^\infty (x+\delta
)^{p+2}e^{-2x}dx,&(6.7)\cr}$$
\noindent where the last step follows from the binomial theorem.\par 
\indent We approximate $\Gamma$ by the Hankel
matrix $\Gamma_N=[\gamma_{j+k-1}{\bf I}_{\{(j,k):j+k\leq N\}}]$, which is zero outside
the top left corner and has rank less than or equal to $N$. By
considering approximation numbers [7], we find that 
$$s_N\leq \Vert \Gamma -\Gamma_{N}\Vert_{B(\ell^2)},\eqno(6.8)$$
\noindent and since the norm of a matrix is smaller than the absolute
sum of its entries, we deduce that 
$$\eqalignno{s_N&\leq \sum_{j=N+1}^\infty j\vert \gamma_j\vert\cr
 &\leq
N^{-p}\sum_{j=1}^\infty j^{p+1}\vert \gamma_j\vert\cr
&\leq {{M_1(\delta )}\over {2\delta^{p+2}N^p}}\Bigl( (2\delta
)^{p+2}+\Gamma (p+3)\Bigr),&(6.9)\cr}$$
\noindent where we have used (6.7) at the
last step. We can repeat the argument with $p+1$ instead of $p$, and
hence deduce that $j^ps_j\rightarrow 0$ as $j\rightarrow\infty .$\par
\indent (ii) We introduce the counting function
$n(t)=\sharp\{ j:ts_j^2\geq 1\}$, and observe that by (6.9),
$$n(t)\leq c_pt^{1/p}\qquad (t>0)\eqno(6.10)$$
\noindent for some $c_p>0$, since $s_j^2\leq c_pj^{-p}$ for all $j$ and 
$n(t)=0$ for $0<t<\Vert\Gamma_\phi\Vert^{-2} .$ A standard
summation formula then gives
$$\eqalignno{ \log\det (I+x\Gamma_\phi^2)&=\log\prod_{j=1}^\infty
(1+xs^2_j)\cr
&=x\int_0^\infty {{n(t)dt}\over{t(x+t)}}\cr
&\leq c_px^{1/p}\pi {\hbox{cosec}}\, (\pi /p).&(6.11)\cr}$$
\noindent This gives the asserted bound on the growth of $\log \det (I+x\Gamma^2_\phi
)$, and shows that $\det (I+z\Gamma^2_\phi)$ has order
zero.\par
\rightline{$\square$}\par
\vskip.1in
\indent The following applies to Q(vi), Q(vii) and Q(viii) of
Examples 3.3.\par
\vskip.05in 

\noindent {\bf Corollary 6.2.} {\sl Let $\phi$ be either:\par 
\noindent (i) $e^{-x}L_n^{(1)}(2x);$\par
\noindent (ii) $(x+s)^{-1}e^{-(x+s)}$; or\par
\noindent (iii) $(x+s)^{-1/2}K_\nu (2\sqrt {x+s})$ where $s>0$. \par
\noindent Then the eigenvalues of the Hankel
operator $\Gamma_{\phi}$ on $L^2(0, \infty )$ are 
of rapid decay, so $j^ps_j\rightarrow 0$ as
$j\rightarrow\infty$ for all integers $p$.}\par
\vskip.05in
\noindent {\bf Proof.} Theorem 6.1 applies directly to $\phi$ in
cases (i) and (ii), and
we can adapt the proof of Theorem 6.1 to deal with $\phi$ in case
(iii). Let $z=re^{i\theta }$ and $0<\delta <s$. One can show that there
exists a constant $C$ such that 
$$\vert (z+s)^{1/2}K_\nu (2\sqrt{z+s})e^z\vert\leq
C(s-\delta )^{-1/2}e^{r\cos\theta
-(r\cos\theta )^{1/2}}\qquad (\Re z>-\delta ),\eqno(6.12)$$
\noindent and hence that 
$$\vert \psi (z)e^z\vert \leq C(s-\delta )^{-1/2}e^{2x+\delta -(2x+\delta
)^{1/2}}\eqno(6.13)$$
\noindent when $z$ lies on the circle with centre $x$ and radius
$x+\delta $. Now we can apply Cauchy's estimates as in (6.5), and then follow the proof 
of Theorem 6.1.\par
\rightline{$\square$}\par
\vskip.1in

\indent Let the Fourier transform of $f$ be ${\cal F}f(\xi
)=\int_{-\infty}^\infty e^{-i\xi x}f(x)dx/\sqrt{2\pi}.$\par
\vskip.05in

\noindent {\bf Proposition 6.3.} {\sl Let $\phi$ be an analytic function on the
strip $\{ z=x+iy: \vert y\vert <\sigma\}$ and suppose that there
exist $K$, $\delta >0$ such that $\vert
\phi (z)\vert \leq Ke^{-\delta \vert x\vert}$ for all real $x$ and $\vert
y\vert <\sigma$. Then there exist $C_1$ and $\kappa_2>0$ such that the
singular numbers of the Hankel
operator $\Gamma_\phi$ with kernel $\phi (x+y)$ all satisfy}
$$s_N(\Gamma_\phi )\leq C_1e^{-\kappa_2N^{1/3}}.\eqno(6.14)$$
\vskip.05in
\noindent {\bf Proof.} By shifting the line of integration in the
Fourier integral, one can show that ${\cal F}\phi $ is analytic on the
strip $\{\zeta =\xi +i\eta : \vert \eta\vert <\delta \}.$ Further,
there exist $C_3$ and $\varepsilon>0$ such that $\vert {\cal F}\phi (\xi
+i\eta )\vert \leq C_3e^{-\varepsilon\vert \xi\vert}$ for all $\xi
+i\eta$ on this strip.\par
\indent The Fourier transform of $h_n$ is easy to compute, and by
Plancherel's theorem we have
$$ \gamma_n=\int_0^\infty \phi (x)h_n(x)\, dx
={{i}\over{\sqrt{2\pi}}}\int_{-\infty}^\infty {\cal F}\phi (\xi ){{(\xi
-i)^n}\over{(\xi +i)^{n+1}}} d\xi .\eqno(6.15)$$
\noindent We shift the line of integration up to $\{\zeta =\xi
+i\delta :-\infty <\xi <\infty \}$, and split the range of integration; the middle range
of integration $[-\xi_0 , \xi_0]$ is dominated by the rational factor, where we observe that 
$$\Bigl\vert {{(\xi +i\delta -i)^n}\over{(\xi +i\delta
+i)^{n+1}}}\Bigr\vert\leq \Bigl({{\xi_0^2+(1-\delta )^2}\over
{\xi_0^2+(1+\delta )^2}}\Bigr)^{n/2}\qquad (\vert \xi\vert\leq
\xi_0);\eqno(6.16)$$
\noindent meanwhile, the contributions from $\vert\xi\vert\geq \xi_0$ are
bounded by 
$$\Bigl\vert \int_{-\infty}^{-\xi_0}+\int_{\xi_0}^\infty
{{(\xi+i\delta 
-i)^n}\over{(\xi +i\delta +i)^{n+1}}} {\cal
F}\phi (\xi+i\delta )\,d\xi \Bigr\vert\leq
\int_{-\infty}^{-\xi_0}+\int_{\xi_0}^\infty \vert {\cal F}\phi (\xi
+i\delta )\vert d\xi .\eqno(6.17)$$
\noindent In particular, with $\xi_0=n^{1/3}$, we obtain the bound
$$\vert \gamma_n\vert\leq 2n^{1/3}C_3\Bigl( {{n^{2/3}+(1-\delta
)^2}\over{n^{2/3}+(1+\delta)^2}}\Bigr)^{n/2}+{{2C_3}\over{\varepsilon}}
e^{-\varepsilon n^{1/3}}.\eqno(6.18)$$
\noindent Due to a similar approximation argument as (6.7), the singular numbers of $\Gamma_\phi$ 
consequently satisfy (6.14).\par
\rightline{$\square$}\par
\vskip.05in
\indent This result applies to an orthonormal sequence of functions
which is used to analyse the Gaussian orthogonal ensemble as in
[1].\par
\vskip.05in

\noindent {\bf Corollary 6.4.} {\sl Let $\phi_n$ be the $n^{th}$
Hermite function. Then the singular numbers of the Hankel operator
$\Gamma_{\phi_n}$ on $L^2(0, \infty )$ with kernel $\phi_n(x+y)$ are
of rapid decay as in (6.14).}\par
\vskip.05in
\noindent {\bf Proof.} The hypotheses of Proposition 6.3 hold for
$\phi_n (z)=H_n(z)e^{-z^2/4}$, where $H_n$ is the Hermite polynomial of
degree $n$.\par
\rightline{$\square$}\par

\vskip.05in
\noindent {\bf References}\par
\noindent [1] G. Aubrun, A sharp small deviation inequality for the
largest eigenvalue of a random matrix, in: S\'eminaire de
Probabilit\'es XXXVIII, in: Lecture Notes in Math.,
vol. 1857, Springer, Berlin, 2005, pp. 320--337.\par
\noindent [2] G. Blower, {Operators associated with the soft and hard edges from unitary
ensembles}, {J. Math. Anal. Appl.} {337} (2008) 239--265. doi:10.1016/j.jmaa.2007.03.084.
\par
\noindent [3] G. Blower, {Integrable operators and the squares
of Hankel operators, {J. Math. Anal. Appl.} j.jmaa2007.09.034.\par
\noindent [4] I.S. Gradshteyn and I.M. Ryzhik, Table of Integrals,
Series, and Products, Fifth Edition, Academic Press, San Diego, 1965.\par
\noindent [5] F.W.J. Olver, Asymptotics and Special Functions,
Academic Press, New York, 1974.\par
\noindent [6] V.V. Peller, Hankel Operators and Their Applications,
Springer, Berlin, 2004.\par
\noindent [7] A. Pietsch, Eigenvalues and $s$-Numbers, Cambridge 
University Press, Cambridge, 1987.\par
\noindent [8] E.C. Titchmarsh, Eigenfunction Expansions Associated with
Second Order Differential Equations I, Clarendon Press, Oxford,
1962.\par 
\noindent [9] C.A. Tracy and H. Widom, {Level-spacing distributions and the Airy
kernel,} {Comm. Math. Phys.} {159} (1994) 151--174.\par
\noindent [10] C.A. Tracy and H. Widom, {Level spacing distributions and the Bessel
kernel,} {Comm. Math. Phys.} {161} (1994) 289--309.\par
\noindent [11] C.A. Tracy and H. Widom, {Fredholm determinants, differential
equations and matrix models}, {Comm. Math. Phys.} {163} (1994) 33--72.\par
\noindent [12] E.T. Whittaker and G.N. Watson, A Course of Modern
Analysis, Fourth edition, Cambridge University Press, Cambridge, 1965.\par

\vfill
\eject
\end